\documentclass{article}
\usepackage{amsmath,amsfonts,amsthm,amssymb,amscd,colordvi}

\newcommand{\e}{\varepsilon}
\newcommand{\vk}{\varkappa}

\newcommand{\vp}{\varphi}

\newcommand{\R}{{\mathbb R}}
\newcommand{\C}{{\mathbb C}}

\newcommand{\Z}{{\mathbb Z}}

\newcommand{\T}{{\mathbb T}}

\newcommand{\XX}{{{X_m^T}}}
\newcommand{\YY}{{{Y_m^T}}}
\newcommand{\ZZ}{{{Z_m^T}}}
\newcommand{\HH}{{{\cal H}^m}}
\newcommand{\tb}{{\widetilde\Box}}
\newcommand{\fu}{{f(t,x,u,\nabla u, \dot u)}}
\newcommand{\Om}{{\cal O}}

\newcommand{\Omo}{{\cal O}^0}

\newcommand{\Pt}{{}_\theta\!\Phi}
\newcommand{\Omdo}{ {}^\delta\! \Om^0  }
\newcommand{\Sth}{{}_\theta\!S }
\newcommand{\Sta}{{}^\tau\!S}
\newcommand{\Sde}{{}^\delta\!S}
\newcommand{\Oth}{{}_\theta {\cal O}}
\newcommand{\Tt}{{}_\theta\! R}
\newcommand{\Rd}{{}^\delta\! R}
\newcommand{\TT}{{\bf T\!r}\,}
\newcommand{\nbh}{neighbourhood}

\newcommand{\cA}{{\cal A}}

\newcommand{\cH}{{\cal H}}
\newcommand{\cL}{{\cal L}}

\newcommand{\iso}{\,\, \,\overrightarrow{\sim}\, \,\,}

\def\12{\tfrac12}

\theoremstyle{plain}
\newtheorem{theorem}{Theorem}[section]
\newtheorem{lemma}[theorem]{Lemma}
\newtheorem{proposition}[theorem]{Proposition}
\newtheorem{corollary}[theorem]{Corollary}
\theoremstyle{definition}

\theoremstyle{remark}

\newtheorem{example}[theorem]{Example}

\numberwithin{equation}{section}

\begin{document}

\title{Analyticity of solutions for quasilinear wave equations and other quasilinear systems
}

\author{Sergei
Kuksin\footnote{
Universit\'e Paris 7
UFR de MathŽmatiques, 
75205 Paris Cedex 13, France, e-mail:
kuksin@gmail.com }, Nikolai Nadirashvili\footnote{LATP, CMI,
39, rue F. Joliot-Curie, 13453 Marseille
France,
e-mail: nicolas@cmi.univ-mrs.fr}}
 \maketitle

\date{}
%%\date{(preliminary version)}
\maketitle

\begin{abstract}
We prove the persistence of analyticity for classical solution of the  Cauchy problem
for quasilinear wave  equations  with analytic data. Our results show that the analyticity of solutions,
stated by the  Cauchy-Kowalewski and  Ovsiannikov-Nirenberg theorems, lasts till a classical 
solution exists. Moreover, they show that if the equation and the Cauchy data are
analytic only in a part of space-variables, then a classical solution also is analytic in these variables. 
The approach applies to other quasilinear equations and implies  the persistence
of  the space-analyticity (and the partial space-analyticity) 
of their classical solutions.  
\end {abstract}
\bigskip

\setcounter{section}{-1}

\section{Introduction}
\label{s0}

Consider a quasilinear wave equation:
\begin{equation}\label{01}
%\Box u+f(u,\nabla u,\dot u)=0 ,
\Box u+f(t,x, u,\nabla u,\dot u)=0 , \quad  \text{dim\,}x=d,\;\; 
t\in\R,
\end{equation}
\begin{equation}\label{02}
u_{t=0}=u_0,  \qquad {\dot u}_{t=0}=u_1, 
\end{equation}
where 
$f$ is a real-analytic function of all its arguments  % in the arguments $u, \nabla u, \dot u$, 
 and the Cauchy data $u_0, u_1$
are real-analytic functions of $x$. To begin with we assume the periodic boundary conditions:
$$
x\in \T^d=  \T^d_\Gamma= \R^d/ \Gamma \qquad (\Gamma\; \text{ is a lattice)}.
$$
Regarding the solvability of the Cauchy problem \eqref{01}, \eqref{02} two facts are well known:
the life-span of its classical solution %to the Cauchy problem %for general nonlinear hyperbolic equation
 is non zero, i.e.,  there is $T>0$  such that on $[0,T)$ the problem \eqref{01}, \eqref{02}
has a $C^2$-solution $u$, e.g. see \cite{Hor} and Proposition~\ref{prop} below. 
On the other hand by the Cauchy-Kowalewski theorem \cite{Kow}
there is a positive ${\e}_1 $ such that for $t\in [0,{\e}_1 )$ the solution $u$ is
real-analytic.  %The latter applies to more general equations. In particular, to quasilinear equations with a nonlinearity, 
%depending on $(t,x)$:
%\begin{equation}\label{01}
%\Box u+f(t,x, u,\nabla u,\dot u)=0 .
%\end{equation}
  Ovsiannikov and Nirenberg gave a beautiful generalization of the latter  % Cauchy-Kowalewski
   theorem
to  equations \eqref{01}, where the nonlinearity $f$ is continuous in $t$ (and still is analytic in other variables), 
see in \cite{Nish}. By their result, for any real-analytic $u_0$ and $u_1$ there is a positive ${\e}_2 $ such that for $t\in [0,{\e}_2 )$ the solution $u$ is real-analytic in $x$. 
  
 From the proofs of the Cauchy-Kowalewski  and   Ovsiannikov-Nireberg theorems all that one can say
 about the life-spans of analyticity  ${\e}_1$ and ${\e}_2$ is just their positivity. 
 However,  for the different classes of quasilinear wave
 equations, supplemented with sufficiently smooth initial data,
  the life-span $T$ of classical (smooth) solutions  is often fairly large,
 sometimes $T=\infty $. The  natural question  is if the range of analyticity is extendable up
 to $T$. 
 
 By a general result of Alinhac and Metivier \cite{AM84} 
  the life-span of analyticity ${\e}_1 $ in the Cauchy-Kowalewski theorem
is equal to $T$. The proof of their  theorem is very technical and  involves complicated  recombination of the Taylor's coefficients.  
 (In \cite{BB}  a similar result was obtained earlier for solutions of the 2d Euler equation, using  hyperbolic features of that equation.) 
The aim of this paper is to give a short and transparent proof  of this (and actually more general) properties of 
solutions of quasilinear wave equations.  We also show  that the life-span of analyticity ${\e}_2 $ in  the Ovsiannikov-Nirenberg theorem equals  the life-span time $T$.

Let $f$ be continuous in $t$, $H^m$-smooth in $x$ and analytic in $u,\nabla u, \dot u$ (see Section~\ref{s11}). 
Let $u(t,x), \ 0\le t\le T, x\in\T^d$,  be a solution of the Cauchy problem  \eqref{01}, \eqref{02} which defines a
continuous curve $t\mapsto u(t,\cdot)\in H^{m+1}(\T^d)$,  $\ C^1$-smooth as a curve in $H^m(\T^d)$. Assume that 
 $m>d/2$  (so $u(t,x)$ is a continuos function, $C^1$-smooth in $t$).

\begin{theorem}\label{t01}
 {\it %Let $u(t,x)$, where $0\le t\le T, x\in\T^d$,   be a solution of the Cauchy problem
 %\eqref{01}, \eqref{02}, $H^{m+1}$-smooth in $x$ and $C^1$-smooth in $t$.  Then 
 
 i) If $f$ and $u_0, u_1$ are real-analytic in $(x_1,\dots,x_k)$, $1\le k\le d$, then $u$ also is  real-analytic in 
 these variables,
 
 ii) if $f$ and $u_0, u_1$ are   real-analytic in all their arguments, then $u$ also is.}
 \end{theorem}
 \smallskip
 
 Note that the first assertion of the theorem and the local in time existence of a classical solution imply that if $u_0,u_1$ and $f$ 
 are sufficiently smooth in $x$, continuous in $t$ and analytic in $x_1,\dots,x_k,  u, \nabla u$ and $\dot u$, then the problem 
  \eqref{01}, \eqref{02}  has a unique local in time solution,  analytic in $x_1,\dots,x_k$ (see below Corollary~\ref{c2}).
This generalises  the Ovsiannikov-Nirenberg theorem for equations of such class. 
%is a   generalisation of  

 Theorem~\ref{t01} is proved in Section 1; its   proof 
   is based on properties of the   nonlinear semigroup, generated  by the problem  \eqref{01}, \eqref{02}. 
 In  Theorem~\ref{t3}, Section~2,  we show  that the assertion  holds for  solutions of  \eqref{01}, \eqref{02} defined locally, in a characteristic  cone in $\R\times\R^{d}$. The local result on  the analyticity  implies the analyticity of 
 global solutions defined on the whole torus.  
 It  straightforwardly generalizes to equations on homogeneous spaces and implies  the corresponding global results. For example, Theorem~\ref{t01}.ii) remain true for equations in the   standard sphere $S^d$, see Section~\ref{s31}.

 We preface Theorem~\ref{t01}  to  local Theorem \ref{t3}
 since the assertion~i) of  the former and its proof remain true for other classes of equations for
  which the latter  is no more  valid. E.g. see Section~\ref{s3} for 
 quasilinear parabolic equations, the 3d~Navier-Stokes system and NLS equations. 
 In the same time,  the proof of assertion~ii) does not generalise to 
  non-hyperbolic equations (and for quasilinear parabolic equations its claim is wrong).  So while the
   Cauchy-Kowalewski theorem is an   assertion about hyperbolic equations, 
  the Ovsiannikov-Nirenberg theorem   describes a general property of a large class of quasilinear systems.

  We note that similar $C^{\infty }$-smooth  properties
of solutions for  \eqref{01}, \eqref{02} are known, see  \cite{Koch, Sog}.

 Our proof relies heavily on the   theory of analytic mappings in Banach spaces.
  For the reader's convenience we summarised  it in an   appendix to this paper.

 {\bf Acknowledgments.} The authors would like to thank H. Koch for very  useful discussions. SK was
 supported   l'Agence Nationale de la Recherche through the grant  ANR-10-BLAN 0102.

\section{Global results: quasilinear wave equation on $\T^d$} \label{s1}
\subsection{Single equation.} \label{s11}
Here we  study  the Cauchy problem for 
a quasilinear wave equation \eqref{01}, \eqref{02} on $\T^d=\R^d/\Gamma$, 
where the function $f$ is continuous in all variables, is 
 $H^m$ smooth in $x$, where 
$$
m>d/2\qquad \text{(or  $m>d/2-1$ if $f$ does not depend on $\nabla u$ and $\dot u$)}, 
$$
and is (real-)analytic in the arguments $u, \nabla u, \dot u$. \footnote{More precisely, the function  $f(t,x,u,\xi,\eta)$ is 
continuous and analytically in ${\frak z}=(u,\xi,\eta)$ extends to a \nbh  \ $Q$ of $\R^{d+2}$ in $\C^{d+2}$. The 
extended function is bounded uniformly on sets $[0,T]\times \T^d\times Q_R$, where 
$Q_R=Q\cap\{{\frak z}\in\C^{d+2}: |{\frak z}|<R\}$.  Moreover, $\|f(t,\cdot, {\frak z})\|_m\le C(R)$ for all 
$t\in[0,T]$ and ${\frak z} \in Q_R$, for each $R>0$.}
We denote by $H^m$ the Sobolev spaces $H^m(\T^d)$ with the norm
$\ 
\|u\|_m=\big(|\nabla^m u|^2_{L_2}+|u|^2_{L_2})^{1/2},
$
and abbreviate $H^{m+1}\times H^m=\HH$. Consider the Cauchy operator for the linear wave equation:
\begin{equation}\label{Cau}
\widetilde\Box: u\mapsto (u_{t=0}, {\dot u}_{t=0}, \Box u).
\end{equation}
It is well known that for any reasonable domain of definition this map is an embedding. For any $T>0$ consider the spaces
$$
\XX=C(0,T;H^{m+1})\cap C^1(0,T;H^m)  ,\qquad 
\YY=\HH\times C(0,T; H^m). 
$$
It is also well known (e.g., see \cite{Tem}) that the inverse operator defines a continuous mapping
\begin{equation}\label{3}
{\widetilde\Box}^{-1}:\YY \to  \XX
\end{equation}
(but certainly $\widetilde\Box$ does not map $\XX$ to $\YY$). 

 The spaces $\YY$ and $\XX$ suit well to study
solvability of the problem  \eqref{01}, \eqref{02}. Indeed, denote $\vk=|(u_0,u_1)|_{\HH}$,  assume that $u_0,u_1$ are 
smooth and that $u(t,x)$ is a smooth solution of the problem such that 
$$
|U(t)|_{\HH}\le 2\vk\qquad \forall\,0\le t\le T'
$$
with some $T'>0$, where $U(t)=(u(t),\dot u(t))$. Taking the $H^m$~scalar-product of  \eqref{01} with $\dot u(t)$, 
we get that 
$$
\frac12 \frac{d}{dt}\|\dot u\|^2_m+C\frac{d}{dt}\|u\|^2_{m+1}\le
C_1\|u\|_m\|\dot u\|_m
+C_2\|\fu\|_{m}\|\dot u\|_{m},
$$
for suitable contants $C, C_1, C_2$. 
By the apriori assumption and since the space $H^m$ is an algebra, 
 for $0\le t\le T'$ the r.h.s. is bounded by $C(\vk)$. Therefore $|U(t)|_{\HH}^2\le\vk^2+tC_1(\vk)$.
So there exists $T(\vk)>0$ such that  $|U(t)|_{\HH} \le 2\vk$ for $0\le t\le T(\vk)$. In the usual way the obtained 
apriori estimate implies

\begin{proposition}\label{prop}
There exists $T'>0$, depending only on  $f$ and $|(u_0,u_1)|_{\HH}$, such that the problem  \eqref{01}, \eqref{02} has a 
unique solution $u(t,x),\ 0\le t\le T'$, belonging to the space $Y_m^{T'}$. 
\end{proposition}

It is well known that in general the local solution $u(t)$ cannot be extended to all $t\ge0$.  The construction below 
gives a convenient implicit description of the set of initial data for which a solution exists for $0\le t\le T$. This 
construction is a part of the PhD thesis of the first author \cite{K0}. 

Denote ${\widetilde\Box}^{-1}\YY=\ZZ$ and provide $\ZZ$ with a norm, induces from 
$\YY$ by ${\widetilde\Box}^{-1}$. This is a Banach space, 
\begin{equation}\label{4}
\widetilde\Box:\ZZ\to \YY \qquad \text{is an isomorphism},
\end{equation}
\begin{equation}\label{5}
\ZZ\subset     \XX
\qquad \text{continuously}
\end{equation}
by \eqref{3}, and $X^T_{m+1}\subset \ZZ$. 
 Denote by $F$ the nonlinear  differential  operator 
$F(u)=\fu$
 and  by  $\Phi$ -- the operator of the Cauchy problem
\eqref{01}, \eqref{02}.   That is 
 \begin{equation}\label{oper}
 \Phi(u)=\tb(u)+  (0,0,F(u)).
 \end{equation}
 Since $m>d/2$, then the space $ C(0,T; H^m)$ is a Banach algebra. Using \eqref{4} and \eqref{5}
 we get that the mapping
\begin{equation}\label{6}
\Phi:\ZZ\to \YY
\qquad \text{is analytic}
\end{equation}
(e.g.  see \cite{RS}; 
cf. Exemple \ref{eA1} in  Appendix). 
It is well known  that the Cauchy problem \eqref{01}, \eqref{02} with zero in the r.h.s.
replaced by any function from $ C(0,T; H^m)$  has at most one solution in  $\XX$.  So $\Phi$ is an embedding.
Consider its differential in any point $u\in\ZZ$:
\begin{equation}\label{7}
d\Phi(u)(v) = \big(v_{t=0}, \dot v_{t=0}, \Box v+d_3f[u]v+d_4f[u]\nabla v +d_5f[u]\dot v\big).
\end{equation}
Here $f[u]=f(x,u,\nabla u,\dot u)$ and $d_j$ denotes the differential with respect to the $j$-th variable. 
\begin{lemma}
\label{l1}
For any $u\in\ZZ$ the map $d\Phi(u):\ZZ\to \YY$ is an isomorphism. 
\end{lemma}
\begin{proof}
%It is clear that the mapping is an embedding. 
 For any $(v_0,v_1,g)\in\YY$ consider the corresponding Cauchy problem which we write as
\begin{equation}\label{8}
\Box v+v+\big(d_3f[u]v+d_4f[u]\nabla v +d_5f[u]\dot v -v
\big)=g,\quad v(0)=v_0,\ \dot v(0)=v_1. 
\end{equation}
Let $v_0,v_1$ and  $g$ be smooth and $v$ be a smooth solution of the problem. Multiplying the equation by $\dot v$ in $H^m$ and using that 
the space $H^m$ is an algebra, we get:
$$
\frac12\frac{d}{dt}\|\dot v\|^2_m+C  \frac{d}{dt}\|v\|^2_{m+1}\le
C_2\|v\|^2_m+C_3\|v\|_{m+1}\|v\|_m
+C_4\|\dot v\|_m^2%   \|v\|_m 
 +  \|g\|_m^2, 
$$
where the constants $C_j$ are continuous functions of $\|u\|_{\XX}$. We immediately get from this relation that 
$$
\|v\|_{\XX}\le C(\|u\|_{\XX})\|(v_0,v_1,g)\|_{\YY}.
$$
%In particular, $d\Phi(u)$ is an embedding. 
In the usual way this apriori estimate implies that \eqref{8} has a unique solution $v\in\XX$. Then we see from
\eqref{8} that $\tb v\in\YY$. So $v\in\ZZ$ and  $\|v\|_{\ZZ}\le C'(\|u\|_{\ZZ})\|(v_0,v_1,g)\|_{\YY}$. 
\end{proof}

Since $\Phi$ is an embedding, then Lemma~\ref{l1} jointly with the inverse function theorem (see Theorem~\ref{tA2})
imply
\begin{lemma}
\label{l2}
The mapping $\Phi$ is an analytic diffeomorphism of the space $\ZZ$ and a domain ${\cal O}\in \YY$.
\end{lemma}

Therefore if for some $(u_1,u_2,g)\in\YY$ the problem \eqref{01}, \eqref{02} with zero in the r.h.s. of
\eqref{01} replaced by $g$ has a solution $u\in\XX$, then $(u_1,u_2,g)\in\cal O$, 
$u$ belongs to $\ZZ$ and analytically depends on
$(u_1,u_2,g)$. Denote
\begin{equation}\label{9}
\Omo=\{(u_0,u_1)\in\HH :   (u_0,u_1,0)\in\Om\}.
\end{equation}
Then  for 
$0\le t\le T$ the flow-maps 
$$
S_0^t:\Omo\to\HH, \qquad 
(u_0,u_1)\to (u(t), \dot u(t)),
$$
 are well defined and analytic.

 We recall that $T$ is any positive number and that the domain $\Omo$ depends on the time-interval $[0,T]$,   
 $\Omo=\Omo([0,T])$. Similar we may study solutions of \eqref{01}, \eqref{02} on negative time-intervals $[-T,0]$. 
 The assertions above remain true for  operators $S_0^t$ with $t\in[-T,0]$ and with the domain $\Omo= \Omo([-T,0])$. 
 Finally, we may consider eq.~\eqref{01} with the Cauchy data given 
  not at $t=0$, but at $t=t_1$,  for  arbitrary  $t_1\in[0,T]$. In this way we find
  that the flow-maps $S_{t_1}^{t_2}$, where $t_1, t_2\in[0,T]$, are analytic operators with domains
  $\Om^{t_1} ( [0,T])$.

 It is clear from our construction that the operators $S_{t_1}^{t_2}$ and $S_{t_2}^{t_1}$ with domains  $\Om^{t_1}([0,T])$
 and  $\Om^{t_2}([0,T])$  are inverse to each other, and that an operator  $S_{t_1}^{t_2}$ analytically extends to the
 bigger domain $\Om^{t_1}([t_1,t_2])$. This is its maximal domain of definition.

 \subsection{Families of equations.}\label{s12}
We fix any $k\in\{1,\dots,d\}$ and assume that 
$$
\T^d=\T^k\times \T^{d-k}
$$
(that is, $\Gamma=\Gamma_k\oplus \Gamma_{d-k}$ and $\T^k=\R^k/\Gamma_k$, $\T^{d-k}=\R^{d-k}/\Gamma_{d-k}$). 
We  make the torus $\T^{k}=\{\theta=(\theta_1,\dots,\theta_{k})\}$ to act on $\T^d$ 
by the shifts  $\Tt$,  
$$
\Tt(x)=(x^I +\theta ,x^{II}),
$$
where 
$
x^I=(x_1,\dots,x_k)$ and $ x^{II}=(x_{k+1},\dots,x_d).
$
Then the torus acts on the operators $F$ by shifting  their  coefficients:
$
(\Tt F)(u)=f(t,\Tt x,u,\nabla u, \dot u).
$
Clearly we have 
\begin{equation}\label{x.1}
  (\Box+ \Tt F)(\Tt u)=\Tt \big((\Box+F)(u)\big).
\end{equation}
 The operators of the  shifted  Cauchy problem 
  $\Pt(u)=\tb u+(0,0,(\Tt F)u)$
define a mapping  
\begin{equation}\label{x.11}
 \bar\Phi^1  :\T^{k}\times \ZZ\to \YY,\qquad (\theta,u)\to \Pt(u). 
\end{equation}

\begin{lemma}
\label{l14}
Assume that the function $\fu$ is analytic in $x^{I}$. Then the mapping 
$\bar\Phi^1$ is analytic. 
\end{lemma}
\begin{proof}
By \eqref{6} and Corollary \ref{cA2}
we only have to check that the mapping is locally bounded and is
 analytic in $\theta$. The local boundedness of this map (and its continuity) follow from the 
 Banach algebra property of the spaces $H^s$ with $s>d/2$, e.g. see \cite{RS}; cf. \eqref{Sob}. 
  Since $f$ is analytic in $x^{I}$, then 
$ \Pt(u)$  complex-analytically depends on $\theta$ from the complex vicinity of $\T^{k}$.
This implies the assertion. 
\end{proof}

By the results of Section \ref{s12}, for any $\theta$ the operator $\bar\Phi^1(\theta,\cdot)$ defines an analytic diffeomorphism 
of $\ZZ$ and a domain $\Oth\subset\YY$. By the implicit function theorem
\begin{equation}\label{x.2}
\text{ the mapping }\;\;
 \Oth\ni\xi\to \big(\bar\Phi^1(\theta,\cdot)\big)^{-1}\in\ZZ\quad\text{is analytic in $\xi$ and $\theta$.}
\end{equation}
Denoting $\Oth^0=\{(u_0,u_1) : (u_0,u_1,0)\in\Oth\}$ we see that 
 for $\theta\in\T^{k}$  and $0\le t\le T$ the time-$t$ 
flow-mapping,   corresponding to the nonlinearity $\Tt F$, 
is 
\begin{equation}\label{x.3}
\text{an analytic  transformation}\;\;
 \Sth^t_0
 :\Oth^0\to \HH \;\; \text{which analytically depends on} \  \ \theta.
\end{equation}
Relation \eqref{x.1}, where $(\Box+F)u=0$, implies that 
\begin{equation}\label{x.22}
 \Sth_0^t
 \circ \Tt=\Tt\circ S_0^t.
\end{equation}
In particular, $\Tt{\cal O}^0= \Oth^0$. 
\medskip

Similar for any $\delta\in\R$ we define
$$
(\Rd F)(u)=f(t+\delta,x,u,\nabla u,\dot u), \qquad
{}^\delta\!\Phi=\tb+(0,0,\Rd F).
$$
Assume that there exists $\rho>0$ such that for each value of $(x,u,\nabla u,\dot u)$ the function $t\mapsto\fu$ analytically 
extends to the segment $[-\rho, T+\rho]$. Then the mapping 
$$
\bar \Phi^2:(-\rho, \rho)\times \ZZ\to\YY, \qquad
(\delta, u)\mapsto {}^\delta\!\Phi
(u),
$$
is analytic, and for any $|\delta|<\rho$ it defines an analytic isomorphism 
$$
{}^\delta\!\Phi
 : \ZZ\to
 {}^\delta\! \Om  \subset\YY\,,
$$
which analytically depends on $\delta$. We set
$\ 
\Omdo =\{(u_0,u_1)\in\HH:  (u_0,u_1,0)\in    {}^\delta\! \Om  \}
$
and denote by $\Sde_0^t$ the mapping $S_0^t$, corresponding to the operator $\Rd F$. Then 
\begin{equation}\label{x.5}
\text{the operator }\;\;
  \Sde_0^t :\Omdo\to\HH\quad\text{is analytic and analytically depends on} \ \delta\in(-\rho,\rho). 
\end{equation}

\subsection{Analyticity of solutions} \label{s13}
There is a delicate difference between the smoothness (or analyticity) of solutions for a nonlinear
wave equation in time and in space. For instance, there is a number of results which imply for 
a solution a high smoothness in $x$ and only a limited smoothness in $t$, see \cite{Hor}. Another 
example is given by the Ovsiannikov-Nirenberg theorem.  Accordingly below we consider
smoothness of solutions for \eqref{01}, \eqref{02}  in $x$ and in $t$ separately.

{\it Space-analyticity}.
Assume that the function $f$ as above  is analytic in $x^{I}$, as well as the initial data 
$u_0$ and $u_1$. Assume also that the problem \eqref{01}, \eqref{02} 
has a solution $u\in\YY$. Then $(u_0,u_1)\in\Oth^0$  and 
by \eqref{x.3} $\Sth _0^t u_0$ with $0\le t\le T$  is well defined for $\theta$
from a small ball $B_\e=\{  |\theta|<\e\}$ and is analytic in $\theta$.
Using \eqref{x.22} we have 
$$
u(t,x^I +\theta,x^{II})=(\Tt\circ S_0^t)(u_0,u_1)(x)= (\Sth_0^t\circ \Tt)(u_0,u_1)(x).
$$
Since $u_0$ is analytic in $x^{I}$, then the mapping 
$\ 
\T^{k}\to \HH,\;\ 
B_\e\ni\theta\mapsto \Tt(u_0,u_1),
$
is analytic. Using \eqref{x.3} we get

\begin{theorem}\label{t1}
Assume that the nonlinearity $f$ and the initial data 
$u_0,u_1$ are analytic in $x^{I}\in\T^{k}$ and the problem \eqref{01}, \eqref{02} 
has a solution $u(t,x)\in\YY$, $0\le t\le T, x\in \T^d$. Then $u$ is analytic in $x^{I}$.
\end{theorem}

{\it Time-analyticity}. Now assume that  for a suitable $\rho>0$ 
the  function $\fu$, where $-\rho\le t\le T+\rho$,  is analytic in all its arguments, that the Cauchy data 
$u_0(x)$ and $u_1(x)$ are analytic and that the  problem \eqref{01}, \eqref{02}  has a solution $u\in\YY$.
Denote $U(t)=(u(t),\dot u(t))$. By the Cauchy-Kowalewski theorem, the function $u(t,x)$ is analytic for 
$|t|<\e$ and $x\in\T^d$  with  a suitable $\e>0$.  Therefore the curve $ [0,T]\to\HH,\ t\mapsto U(t)$, also 
 is analytic 
for $|t|<\e$. For any $t_*\in[0,T]$ we write the solution $U(t)$ for $t$ close to $t_*$ as 
$\ 
U(t_*+\tau)=\Sta_0^{t_*}
\circ U(\tau).
$
 Using \eqref{x.5} we get
\begin{lemma}\
\label{l4}
Under the above assumptions   the curve
$
[0,T]\to\HH,\;\  t\to (u,\dot u)(t),
$
is analytic. In particular, $u(t,x)$ is analytic in $t\in[0,T]$ for each $x$. 
\end{lemma}

This result and  Theorem~\ref{t1}  imply 
\begin{theorem}\label{t2}
Assume that the nonlinearity  $\fu$, where $-\rho\le t\le T+\rho$, 
 and  the initial data $u_0,u_1$ are analytic in all variables
and the problem \eqref{01}, \eqref{02} 
has a solution $u(t,x)\in\YY$. Then $u(t,x)$ is an analytic function  in all its arguments.
\end{theorem}

We recall that for any $t_1, t_2\in[0,T]$ the flow-map $S_{t_1}^{t_2}$ defines an  analytical isomorphism  
$\Om^{t_1}([0,T])  \iso \Om^{t_2}([0,T]) $. Denote by  $\cA(\T^d)$ the space of analytic functions on $\T^d$.

\begin{corollary}\label{c1}
If the nonlinearity $f$ is analytic in $x, u, \nabla u, \dot u$, then the 
 mapping $S_{t_1}^{t_2}$  defines a bijection $\Om^{t_1}([0,T]) \cap\cA(\T^d) \iso  \Om^{t_2}([0,T])\cap\cA(\T^d) $.
 If $f$ is also analytic in $t\in[-\rho, T+\rho]$, then  for each $u\in \Om^{t_1}([0,T])\cap\cA(\T^d) $ the curve  $(t_1,t_2)\mapsto S_{t_1}^{t_2}(u)$ is analytic in $t_1$ and $t_2$. 
\end{corollary}

Combining Theorem~\ref{t1} with Proposition~\ref{prop} we get 

\begin{corollary}\label{c2}
Let the nonlinearity $f$ and the Cauchy data $u_0,u_1$ be as in Theorem~\ref{t1}. Then there exists $T'>0$ such that 
for $0\le t\le T'$ the Cauchy problem has a unique solution $u\in Y_m^{T'}$, which is analytic in $x^I$. 
\end{corollary}

The global results above generalise to other classes of quasilinear PDE. E.g., to 
quasilinear parabolic and Schr\"odinger equations, see Section~\ref{s3}.  
%Moreover, they remain true for strongly nonlinear hyperbolic equations.
%This will be shown in a separate publication 

\section{Local results}\label{s2}
\subsection{Equations in characteristic cones.}
In this section we consider the problem  \eqref{01}, \eqref{02} defined in a characteristic cone
in $\R^{d+1}$. 

 Let $T>0$ and $ 0<a<T$.
Denote by $K$ a truncated characteristic cone:
\begin{equation*} %\label{cone}
K=K(T,a)=\{ (t,x)\in [0, T-a]\times \R^d: |x|\leq T-t \}
\end{equation*}
(below for short we call it  a cone).
Denote by $B_r\subset \R^d$ the closed  ball of radius $r$ centered in the origin and denote 
%Let $f\in H^s(B_2$ and
 $b^t=B_{T-t}$. In this section we study the problem \eqref{01}, \eqref{02} in the cone $K$, where the Cauchy data 
 are given on the ball $K\cap \{t=0\}$, identified with $b^0$. The nonlinearity $f$ is assumed to be analytic in all 
 its variables and analytically extendable  to $U_\e(K)\times \R^{d+2}$, where  $U_\e(K)$ is the $\e$-vicinity of
 $K$ in $\R^{d+1}$, $\e>0$. 
  Then for  a given Cauchy data the problem \eqref{01}, \eqref{02} has at most
one classical solution, \cite{Hor}. 
 Our goal is to prove for this solution  Theorem~\ref{t01},    assuming that the Cauchy data also are analytic.

Denote by $T\!r_\rho(g)$ the restriction of a function $g(x)$ to the ball $B_\rho$, $\rho>0$.  It is well known  that  for 
%integer
  $k\ge 0$ there exists  a  bounded linear map 
$$L^k_1: H^k(B_1)\rightarrow H^k_0(B_2)$$
such that $T\!r_1\circ L^k_1=\,$id  and $L^k_1u(x)=0$ for $|x|\ge \tfrac32$.
%   and each $u$, where $r<2$ depends only on $L_1$.
 For $0<\rho\le T$ we denote by  $L^k_\rho$ 
the  linear operator
$
L^k_\rho :H^k(B_\rho)\rightarrow H^k_0(B_{2\rho})\subset   H^k_0(B_{2T}),
$
obtained from $L_1^k$ by the dilation (so  $T\!r_\rho\circ L^k_\rho=\,$id).  
 Denote by $C ([0,T-a]; H^{k}(b^t))$ the space of functions $u(t,x)$ on the cone $K$, 
 satisfying\footnote{This space is formed by 
 restrictions to $K$ of functions from $C ([0,T-a]; H^{k}(B_{2T}))$.}
 $$
 \cL^k u \in C([0,T-a]; H_0^{k}(B_{2T}))\quad \text{where} \quad \big( \cL^k u) (t) =L^k_{T-t}(u(t)),
 $$
 and provide this space with the norm, induced  from $ C ([0,T-a]; H^{k}(B_{2T}))$.  Next, for
  $m>d/2$ denote
$
\cH^m=H^{m+1}(b^0)\times H^m(b^0)
$
and set 
\begin{equation*}
\begin{split}
&X'_m= \{u\in C([0,T-a]; H^{m+1}(b^t)) :   \dot u\in C ([0,T-a]; H^m(b^t))\},\\
&Y'_m=\cH^m \times C([0,T-a]; H^m(b^t)).
\end{split} 
\end{equation*}
We denote by $\tb'$ the operator of the Cauchy problem for $\Box$ in $K$ 
 which sends any function $u(t,x)\in X'_m$ to
$\tb' u=(u(0), \dot u(0), \Box u) \in Y'_{m-1}$.

Consider the torus $T^d=\R^d/(4T)\Z^d$,  the corresponding spaces $\XX, \YY, \ZZ$ and the operator $\tb$. 
 We will identify functions on $T^d$ with $4T$-periodic functions on $\R^d$. Denoting by $\TT$  the operator of restricting a function on $[0,T-a] \times \R^d$ to $K$, we get the mappings 
$$
\TT:\XX \to X'_m,\quad \TT:\YY\to Y'_m.
$$

For $\rho\le T$ and any $u(x)\in H_0^s(B_{2\rho})$ denote by 
$\iota$ the operator which first extends $u(x)$ to the cube $[-2T,2T]^d$  by zero  outside the ball $B_{2\rho}$ 
and next extends it
to a $4T$-periodic function. This is a bounded linear operator from $H_0^s(B_{2\rho})$ to  $H^s(T^d)$ for any $s\ge0$ and  any  $0<\rho\le T$.   For a function $u(t,x)$ we set
$
(\iota u)(t,x)=\iota(u(t,\cdot))(x),
$
if the r.h.s. is defined.  Clearly $\TT\circ\iota=\,$id on the space $C([0,T-a]; H^m(b^t))$.

\begin{lemma}\label{l31}
The inverse operator $(\widetilde\Box')^{-1}$  equals
\begin{equation}\label{relation}
(\tb')^{-1}=\TT\circ \tb^{-1}\circ\iota. 
\end{equation}
It  defines a continuous mapping
$(\widetilde\Box')^{-1}:Y'_m \to  X'_m.$
\end{lemma}
\noindent
{\it Proof.} For $(u_0, u_1, g)\in Y'_m$
 consider functions $\hat u_j(x)=\iota\big(L^m_T(u_j(x))\big), j=0,1$ and $\hat g(t,x)= \iota\big( \cL^m(g)\big)$.
 %and extend them, as functions of $x$, byzero to the whole $R^d$. 
 Then $|(\hat u_0,\hat u_1, \hat g   )|_{\YY } \le C |(u_0,u_1,g)|_{Y'_m}$. The solution of the Cauchy 
 problem for $\Box$ 
 in  $T^d$,  $U(t,x)=\tb^{-1}(\hat u_0, \hat u_1, \hat g)$, satisfies 
 $\ 
 |U|_{\XX} \le C |(\hat u_0, \hat u_1, \hat g)|_{\YY }.
 $
% where the spaces $\XX)$ and $Y'_m({T^d})$ are obtained from $X'_m$ and $Y'_m$ by replacing the
 %balls $b^t$ with the whole $R^d$ (cf. Section~\ref{s1}). 
Since solutions of the wave equation in $K$  depend only on the data in the 
characteristic cone, then $(\tb')^{-1}(u_0, u_1, g)$ equals to the restriction of $U$ to $K$. This implies the assertions. 
\qed 
\smallskip

As above, we define a Banach space $Z'_m$ as $(\widetilde\Box')^{-1}Y'_m$,
$\,Z'_m\subset X'_m$. Due to \eqref{relation}  $Z'_m= \TT(\ZZ)$. Denote by $\Phi'$ the operator of the Cauchy problem 
\eqref{01}, \eqref{02} on $K$ (cf. \eqref{oper}). Then   the following diagram is commutative 
\[
\begin{CD}
Z_m @>\Phi>> Y_m \\
@V{\TT}VV          @V{}VV \\
Z'_m @>\Phi'>> Y'_m
\end{CD}
\]
where the second vertical line stands for the mapping $ T\!r_T\times T\!r_T\times\TT $. 
So we derive from  Lemma~\ref{l1} that  for each $u\in Z_m$ the mapping
 $d\Phi' (u)$ is an isomorphism. Hence, $\Phi' $ is an analytic isomorphism
 $$\Phi': Z'_m \iso O,$$ 
 where $O=\Phi'(Z'_m)$ is a domain in $Y'_m$.
   We define  $\Om^0$ by the relation \eqref{9}.  This is a domain 
 in  $\cH^m$ such that the problem
  \eqref{01}, \eqref{02} 
  has a solution $u\in X'_m$ if and only if  
  $(u_0,u_1)\in\Om^0$.  For $0\le \tau\le T-a$ the  mapping
  $$
  S_0^\tau(u_0,u_1)\mapsto (u(\tau), \dot u(\tau)),\qquad \Om^0\to \cH^m,
  $$
 is analytic since $\Phi'^{-1}$ is an analytic mapping on $\Om$. 
 
 Similar we may consider eq.~\eqref{01} on the smaller cone
 $$
K^\tau=\{ (t,x)\in [\tau, T-a]\times \R^d: |x|\leq T-t \}=K\cap ([\tau, T-a]\times\R^d), \quad 0\le\tau<T-a.
$$
In this way we get a domain $\Om^\tau\subset\cH^m(K^\tau)$ such that eq.~\eqref{01} has a solution
$u(t,x)$, $(t,x)\in K^\tau$, which is a trace on $K^\tau$ of a function from $X'_m$, if and only if
$(u(\tau), \dot u(\tau))\in \Om^\tau$.  Clearly the flow-map $S_0^\tau$ is an analytic operator 
$S_0^\tau:\Om^0\to \Om^\tau$. In difference with Section~\ref{s1} this is not an embedding. 

Now we define families of the Cauchy problems and of the corresponding operators $\Phi'$. Since the 
function $f$ analytically in $(t,x)$  extends to $U_\e(K)$, then eq.~\eqref{01} analytically extends to a bigger cone
$$
K^+=\{ (t,x)\in [-\e, T+\e-a]\times \R^d: |x|\leq T+\e-t \},\;\; \e>0.
$$
For any $\theta\in B_\e=\{|\theta|\le  \e\}\subset \R^d $, as before,  we denote by $\Tt$ the shift
$
\Tt(t,x)=(t,x+\theta),
$
 and set
 $$
 K^\theta=\Tt(K)\subset K^+, 
 \qquad (\Tt f)(t,x,u\nabla u, \dot u)=f (t,\Tt x,u\nabla u, \dot u).
 $$
 As before, $\Pt'$ is the operator of the Cauchy problem with the nonlinearity $\Tt f$. The mapping
 $$
 \bar\Phi^1: B_\e\times Z'_m\to Y'_m,\qquad (\theta, u)\to \Pt' (u)
 $$
 is analytic. For each $\theta\in B_\e$ it defines an analytic diffeomorphism 
 $$
  \bar\Phi^1(\theta, \cdot):Z'_m\to \Oth\subset Y'_m,
 $$
 which  analytically depends on $\theta$, as well as its inverse. 
 
 Considering eq. \eqref{01} with the shifted nonlinearity $\Tt f$ in a smaller cone $K^\tau$ we define the corresponding
 domain $\Oth^\tau\subset\cH^m(b_\tau)$, formed by the initial data $(u(\tau), \dot u(\tau))$ for which the shifted equation  has a solution in  $K^\tau$, extendable to a function from $X'_m$. Then the flow-map of
 the shifted equation $\Sth_0^\tau$ is an analytical mapping 
 $$
 \Sth_0^\tau: \Oth^0\to \Oth^\tau
 $$
 which analytically depends on $\theta$, and 
 $\ 
  \Sth_0^\tau\circ \Tt= \Tt\circ S_0^\tau.
 $
 That is, denoting by $u(t,x)$ a solution of \eqref{01}, \eqref{02} we have 
 $$
 u(t,x+\theta)=( \Sth_0^t)(u_0(x+\theta), u_1(x+\theta)). 
 $$
 In particular, if the functions $u_0$ and $u_1$ are analytic in the closed ball $b_0$
  (i.e. analytically extend to its vicinity 
 in $\R^d$), then $u(t,x)$ is analytic in $x$. 
 
 By the Cauchy-Kowalewski theorem the solution $u(t,x)$ is analytic in the vicinity of the disc $b_0\times\{0\}$ 
 in $\R^{d+1}$. Considering shifts of the nonlinearity $f$ by the time-translations and arguing as above (cf.
 Section~\ref{s13}) we find that $u$ is analytic in $t$. We have proven

 \begin{theorem}\label{t3}
 Assume that $f$ is an analytic function on $U_\e (K)\times \R^{d+2}$ for some $\e>0$,  
 and that the Cauchy data  $u_0, u_1$ are analytic
 in $b^0$.  Let the Cauchy problem  \eqref{01}, \eqref{02} have  a solution $u(t,x) \in X'_m$.
 Then $u$ is analytic.
 \end{theorem}
 
 An obvious local version of Theorem~\ref{t1} also is true.
 
 Since the open non-truncated  characteristic  cone
$K^o= 
\{ (t,x)\in [0, T) \times \R^d: |x| < T-t \}
$
is the union of the closed truncated cones $K(T',a)$ with $T'<T$ and $a>0$, then a natural version of Theorem~\ref{t3} holds 
for the cone  $K^o$.

\subsection
{Global problems} \label{s2x}
In the assumptions of Theorems~\ref{t2} let $R>0$ be such that any two points of the ball $B_{R}$ are not
equivalent modulo the lattice $\Gamma$ (which defines the torus $\T^d$). Define $T'=\min\{2R, 2T\}$ and cover  the layer $[0,T']\times \T^d$ by finitely many truncated cones $K(2T', T')$, shifted by
vectors $(0,\xi),\ \xi \in \R^d$.   Then by 
Theorem~\ref{t3} the solution $u$ is analytic in the layer. Iterating this construction (if $T'<T$) we see that 
$u$ is analytic in $[0,T]\times\T^d$. So the  local results of  this section provide  another proof of Theorems~\ref{t1}
and \ref{t2}.

They  also
straightforwardly generalise to  quasilinear wave equations in a connected open domain in
an analytic Riemann  homogeneous space, not necessarily compact.   In this case
$\Box = \partial^2/\partial t^2 - \Delta $, where  $\Delta $ is the corresponding 
Laplace-Beltrami operator. Now the straight  cone $K$
should be replace by the characteristic cone, constructed in terms of the geodesics of the corresponding  metric, and the 
translations $\Tt $ -- by the local isometies. 
  This generalisation implies that  Theorem~\ref{t1}  with $k=d$ 
  and  Theorem~\ref{t2}  remains  true for quasilinear wave 
equations on a compact homogeneous analytic  Riemann  manifold $M$. For example, on the standard sphere $S^d$.

\section{Related results}\label{s3}

\subsection
{Quasilinear parabolic equations} \label{s31}
The approach to study analyticity and partial analyticity of solutions in the space-variables, 
 explained above, applies to other 
equations (to which the Cauchy-Kowalewski and  Ovsiannikov-Nirenberg theorems do not apply). 
For example, to quasilinear parabolic equations
\begin{equation}\label{111}
\dot u -\Delta u+f(t,x,u,\nabla u)=0,\quad x\in \T^d,\;\; t\ge0, \qquad 
u_{t=0}=u_0,
\end{equation}
where $f$ is sufficiently smooth in $t,x$ and is analytic in $u$ and $\nabla u$. 
As in Section~\ref{s11}, one can find  suitable space $\ZZ$ and $\YY$ such that 
the operator $\Phi$ of the Cauchy problem \eqref{111}  defines an analytic diffeomorphism  between $\ZZ$ and 
a subdomain of the space $\YY$, see \cite{K0, K82}.  In the same way as in Section~\ref{s2}  we prove that if $f$ is analytical in  the  space-variables, then 
classical 
solutions of  \eqref{111}  with analytical initial data are space-analytic. This is a well known result, which holds true for $t>0$ 
without assuming 
analyticity of $u_0(x)$; we do not discuss the corresponding  literature but only 
mention the papers \cite{Ang, Ang90}, where analyticity of solution for nonlinear parabolic equations is 
obtained, based on  ideas, related to those in our work and in \cite{K82}. 
But we also can prove that if $f$ is analytic in  $u, \nabla u$ and in a part of the space-variables, as
well as the function $u_0$, then the solution $u(t,x)$   is analytic in these space variables. This result seems new. Note that 
the assertion of Theorem~\ref{t2} does not hold for the problem \eqref{111}, even when $f=0$, since a solution of the Cauchy 
problem \eqref{111}$\mid_{f=0}$ with analytic $u_0(x)$ may be non-analytic in $t$ when $t=0$. 

The approach applies to the Navier-Stokes system on the $d$-torus with $d=2$ or $d=3$,
 perturbed by a sufficiently smooth force $h(t,x)$, see \cite{K82}. It implies that if the initial data and the force $h$ are analytical in space-variables $x_1,\dots, x_k$, where $1\le k\le d$, then a corresponding strong solution $u(t,x)$ remains analytical in 
 this space-variables till it exists. For $k=d$ this is well known, e.g. see in book \cite{DG}.

 Similar consider the 3d NSE in the thin  layer $M \times (0,\e)=\{(\vp,r\})$, where $M=S^2$ or $M=\T^2$. At the boundary
  $(M\times\{0\} ) \cup
( M\times\{\e\})$ impose the Dirichlet or Navier boundary conditions.
   Let the force and initial data are \\
   i) analytic in $\vp$,\\
   ii) bounded uniformly in $t\ge0$, uniformly in $\e\in(0,1)$.\\
   Due to Raugel-Sell  (see \cite{TZ97} and references therein), if $\e> 0$ is sufficiently small,   then 
   there exists a unique strong solution $u(t,\vp,r),\ t\ge0$.  Our result implies  that this solution is analytic in $\vp$.
   \smallskip

    \subsection{NLS equations    }\label{s32}
    The result of Theorem~\ref{t1} remains true for the nonlinear Schr\"odinger equation 
    $$
    \dot u-i\Delta u +f(t,x, \text{Re\,}u, \text{Im\,}u)=0,\qquad u_{t=0}=u_0, \; \; x\in \T^d,
    $$
    where the complex function $f$ is continuous in $t$,  $H^m$-smooth in $x$ ($m>d/2$) and real analytic in 
    $\text{Re\,}u, \text{Im\,}u$.  The proof of the theorem remains literally the same if we choose 
    $X_m^T=C(0,T;H^m)$, $Y^T_m=H^m\times C(0,T;H^m)$ and proceed as in Section~\ref{s1} (cf. \cite{K0}). 
    As before, we can replace $\T^d$ by any homogeneous Riemann space, analytic and compact.
    \smallskip

 \subsection{Smooth and partially smooth solutions}\label{s33}
Results of Sections~\ref{s1} and \ref{s31} concerning spatial analyticity  and partial spacial analyticity 
of solutions remain true, with the same proof, for their spacial smoothness.
 For example if the nonlinearity $f$ and the initial data $u_0,u_1$ of 
the problem  \eqref{01}, \eqref{02}  are smooth in the variables $x_1,\dots,x_k$, $1\le k\le d$, and the problem has a solution 
$u(t,x)\in\YY$, then  $u$ also  is smooth in  $x_1,\dots,x_k$. 
Similar, if for  the Navier-Stokes system on $\T^3$ 
the initial  data and the force are smooth in some variable $x_l$, then a corresponding strong solution is smooth in $x_l$ till it 
exists. 

\section{Appendix: analytic maps between Banach spaces 
}\label{sA}

Notations for  spaces and operators in this appendix are independent from the rest of the paper. Proofs of the two theorems, 
given below, may be found in \cite{PT}, Appendix~B. 

Let $X$ and $Y$ be complex Banach specs, formed by certain classes of complex functions,  and let $X^R$ and $Y^R$ 
be their subspaces, formed by real-valued functions (in the main part of the text we use a number of spaces like that).  Let $O$ be a  domain in $X$. As in the finite-dimensional case, a mapping $F:O\to Y$  is called (complex) analytic if each
pint $x\in O$ has a \nbh \ $U\subset O$ such that the restriction of $F$ to $U$ may be written as a sum of  series
$$
F(x+u)=\sum_{k=0}^{\infty} F^{(x)}_k(u),\qquad x+u\in U. 
$$
Here $ F^{(x)}_k:X\to Y$ is a bounded $k$-homogeneous mapping. That is, it equals to the restriction to the diagonal 
of some bounded and symmetric complex-linear mapping
$
X^k=X\times\dots\times X\to Y. 
$
Moreover,  $|F^{(x)}_k(u)|_Y\le Cr^k|u|_X^k$ for  suitable positive $r$ and $C$ (so the series above converges uniformly 
if $U$ is sufficiently small). 

Definition of a real-analytic mapping $F: X^R\supset O^R\to Y^R$ is similar. It is easy to see that  a mapping $O^R\to Y^R$
is real-analytic if and only if it extends to a complex-analytic map $O\to Y$, defined on some \nbh \ $O$ of
$O^R$ in $X$. The theorem below gives a very convenient  criterion of analyticity:

 \begin{theorem}\label{tA1} 
 A map $F:O\to Y$, where $O$ is a domain in $X$,  is complex-analytic if and only if it is locally bounded \footnote{
 That is, each point $x\in O$ has a \nbh, where $F$ is uniformly bounded.}
 and weakly analytic. The latter means that for each $x\in O$, $u\in X$ and any $\xi\in Y^*$ there exists $\e>0$
 such that the function 
 \begin{equation}\label{a1}
 \{z\in \C, |z|<\e\} \to \C,\qquad  z\mapsto \xi\big( f(x+\e u)\big),
 \end{equation}
 is analytic. 
  \end{theorem}

 \begin{corollary}\label{cA1}  
 If $F$ is continuous, then it is sufficient to require analyticity of maps \eqref{a1} only for a set of triplets $(x,u,\xi)$ which 
 is dense in $O\times Y\times Y^*$.  
  \end{corollary}
\begin{proof}
For any given triplet  $(x,u,\xi)$ as in Theorem \ref{tA1} 
 let us approximate it by a sequence of admissible triplets
$\{(x,u,\xi)_n, n\ge1\}$.  Then the mapping \eqref{a1}  is a limit in the uniform topology of analytical mappings
 \eqref{a1}${}_n$.
So it is analytic. 
\end{proof}

\begin{example}\label{eA1}
Let  $X=H^s(\T^d;\C)$, $s>d/2$, and let $f(x,u)$ be a continuous function on  $\T^d\times Q$, where $Q$ is a \nbh \ of $\R$
in $\C$, analytic in $u$, $H^s$-smooth in $x$ and satisfying 

%i) for each $R$ it is uniformly bounded on $\T^d\times Q_R$, where  $Q_R=\{z\in Q: |z|<R\}$;

i) $\|f(\cdot,u)\|_s\le C(R)$ for each $u\in Q_R=\{z\in Q: |z|<R\}$
 and for  every $R>0$. 
 
\noindent 
Let $O\subset X$ be the  domain  $O=O(R,\delta) := \{u\in X:\|u\|_s<R, \|  \text{Im}(u)\|_s<\delta\}$,
$\delta<R$. 
Then the  mapping $F: u(x)\mapsto f(x,u(x))$ is defined on $O$  if $\delta$ is sufficiently small in terms or $R$ (and $Q$). 
Since $s>d/2$, then $X$ is a Banach algebra and   $F:O\to Y$ is a continuous mapping.
 In particular $F$ is locally bounded. Moreover, if $O'=O(R-\gamma, \delta-\gamma)$, where $\gamma<\delta$, 
then
 \begin{equation}\label{Sob}
 \|F(u)\|_s\le K(s, d, R, C(R), \delta,\gamma)\qquad \forall\, u\in O'.
 \end{equation}
 For integer $s$ it is  not hard to check this directly, for non-integer $s$ see e.g. \cite{RS}.

Let $u(x), v(x)$ and $\xi(x)$ be smooth complex functions on $\T^d$ such that $u\in O$. Consider the function 
 \begin{equation*}
 z \mapsto \int_{\T^d}  f(x,u(x)+z v(x))\xi(x)\,dx,
 \end{equation*}
 defined in the vicinity of $0\in \C$. The integrand is bounded uniformly in $z,x$ and analytic in $z$.  So the function above is
 analytic, and the mapping $F$ is analytic by Corollary~\ref{cA1}. 
 
 What was said implies that if i)  holds and $f$ is real for real arguments, then this function  defines a real-analytical mapping 
 $F:X^R\to X^R$. 
 Consider the space $Y=C(0,T;X)$. For a similar reason $f$ defines a real-analytical mapping 
 $\ 
 Y^R\to Y^R$,  $u(t,x) \mapsto f(x, u(t,x))$.  Same is true if $f$ continuously depends on $t$. 
\end{example}

 \begin{corollary}\label{cA2}  
 If $O_1, O_2$ are domains in complex Banach spaces $X_1, X_2$, then a mapping 
  \begin{equation}\label{a2}
 F: O_1\times O_2\to T
 \end{equation}
 is analytic if and only if it is locally bounded and 
  for each $x_1\in O_1$ and $x_2\in O_2$ the mappings
  \begin{equation}\label{a3}
 O_1\to Y,\;\; x\mapsto F(x,x_2),\qquad   O_2 \to Y,\;\; x'\mapsto F(x_1,x')\
 \end{equation}
 are analytic.
  \end{corollary}
  \begin{proof}
  In one direction the implication is obvious. Now let \eqref{a2} be a locally bounded 
  mapping such that the maps \eqref{a3}
  are analytic. Take any $(x_1,x_2) \in O_1\times O_2$,
  $(u_1,u_2)\in X_1\times X_2$ and $\xi\in Y^*$. Consider the function 
  $\ 
  (z_1,z_2)\mapsto \xi\big( F(x_1+z_1u_1,x_2+z_2u_2)\big), 
  $
  defined in a suitable bi-disc $\{|z_1|<\e_1\} \times \{|z_2|<\e_2\}$.  It is analytic in $z_1$ for $z_2$ fixed, and vice versa. 
  So this function is analytic by the Hartogs theorem, see \cite{Kr}.  Accordingly $F$ is weakly analytic, so it is analytic.
  \end{proof}

For mappings in Banach spaces the analytic implicit function theorem remains true:

 \begin{theorem}\label{tA2} 
 Let a mapping \eqref{a2}  be analytic, and let for some $(x_1^0,x_2^0)\in O_1\times O_2$ its  differential
 in second variable  $d_2 F(x_1^0,x_2^0): X_2\to Y$ be a linear isomorphism. Denote $F(x_1^0,x_2^0)=y^0$.
 Then in a suitable \nbh \ of $(x_1^0,x_2^0)$ in  $O_1\times O_2$ the equation  $F(x_1,x_2)=y^0$ defines
 $x_2$ as a unique analytic function of $x_1$. 
 
   \end{theorem}

\bibliography{biblio}

\providecommand{\bysame}{\leavevmode\hbox to3em{\hrulefill}\thinspace}
\providecommand{\MR}{\relax\ifhmode\unskip\space\fi MR }
\providecommand{\MRhref}[2]{%
  \href{http://www.ams.org/mathscinet-getitem?mr=#1}{#2}
}
\providecommand{\href}[2]{#2}
\begin{thebibliography}{Ang90b}

\bibitem[AM84]{AM84}
S.~Alihnac and G.~Metivier, \emph{Propagation de l'analyticit\'e des solutions
  d'\'equations hyperboliques non-lin\'eaires}, Invent. Math. \textbf{75}
  (1984), 189--204.

\bibitem[Ang90a]{Ang90}
S.~Angenent, \emph{Nonlinear analytic semiflows}, Proc. Roy. Soc. Edinburgh
  Sect. A \textbf{115} (1990), 91--107.

\bibitem[Ang90b]{Ang}
\bysame, \emph{Parabolic equations for curves on surfaces, {P}art {I}. {C}urves
  with p-integrable curvature}, Ann. of Math. \textbf{132} (1990), 451--483.

\bibitem[BB77]{BB}
C.~Bardos and S.~Benachour, \emph{Domaine d'analyticite des solutions de
  l'\'equation d'{E}uler dans un ouvert de ${R}^n$}, Ann. Scu. Norm. di Pisa
  (4) \textbf{4} (1977), 647--687.

\bibitem[DG95]{DG}
Ch. Doering and J.~Gibbon, \emph{{ Applied Analysis of the Navier-Stokes
  equations}}, Cambridge University Press, Cambridge, 1995.

\bibitem[H{\"o}r97]{Hor}
L.~H{\"o}rmander, \emph{{Lectures on Nonlinear Hyperbolic Differential
  Equations}}, Springer-Verlag, Berlin, 1997.

\bibitem[Koc93]{Koch}
H.~Koch, \emph{Mixed problems for fully nonlinear hyperbolic equations}, Math.
  Z. \textbf{214} (1993), 9--42.

\bibitem[Kow75]{Kow}
S.~Kowalewski, \emph{{Zur Theorie der partiellen Differentialgleichungen}}, J.
  Reine Angew. Math. \textbf{80} (1875), 1--32.

\bibitem[Kra92]{Kr}
S.~G. Krantz, \emph{{Function Theory of Several Complex Variables}}, AMS
  Chelsea Publishing, Providence, Rhode Island, 1992.

\bibitem[Kuk81]{K0}
S.~B. Kuksin, \emph{Diffeomorphisms of functional spaces that correspond to
  quasili\-near equations}, PhD Thesis, Moscow State University (1981).

\bibitem[Kuk82]{K82}
\bysame, \emph{Diffeomorphisms of functional spaces that correspond to
  quasili\-near pa\-ra\-bo\-lic equations}, Math. USSR Sbornik \textbf{117}
  (1982), 359--378.

\bibitem[Nis77]{Nish}
T.~Nishida, \emph{A note on a theorem of {N}irenberg}, J. Differential Geom.
  \textbf{12} (1977), 629--633.

\bibitem[PT87]{PT}
J.~P{\"o}schel and E.~Trubowitz, \emph{{Inverse Spectral Theory}}, Academic
  Press, Boston, 1987.

\bibitem[RS96]{RS}
T.~Runst and W.~Sickel, \emph{{ Sobolev Spaces of Fractional Order, Nemytskij
  Operators, and Nonlinear Partial Differential Equations}}, Walter de Gruyter
  and Co, Berlin, 1996.

\bibitem[Sog08]{Sog}
C.~Sogge, \emph{{Lectures on Non-Linear Wave Equations}}, second ed.,
  International Press, Boston, 2008.

\bibitem[Tem97]{Tem}
R.~Temam, \emph{{Infinite-Dimensional Dynamical Systems in Mechanics and
  Physics}}, Springer-Verlag, Berlin, 1997.

\bibitem[TZ97]{TZ97}
R.~Temam and M.~Ziane, \emph{Navier-{S}tokes equations in thin spherical
  domains}, Contemp. Math., AMS \textbf{209} (1997), 281--314.

\end{thebibliography}
\bibliographystyle{amsalpha}
\end{document}